\documentclass[12pt,reqno]{amsart}
\usepackage{amsmath, amsthm, amssymb, stmaryrd}
\usepackage{hyperref}
\usepackage{enumerate}
\usepackage{url}
\usepackage{color}
\usepackage{comment}

\let\OLDthebibliography\thebibliography
\renewcommand\thebibliography[1]{
  \OLDthebibliography{#1}
  \setlength{\parskip}{0pt}
  \setlength{\itemsep}{0pt plus 0.3ex}
}

\usepackage{mathtools}

\usepackage{multirow, array}
\usepackage{placeins}

\usepackage{caption} 
\advance \topmargin by -\headheight
\advance \topmargin by -\headsep
     
\evensidemargin \oddsidemargin
\newtheorem{theorem}{\bf Theorem}[section]

\newtheorem{lemma}[theorem]{\bf Lemma}

\numberwithin{equation}{section}

\newcommand*\wrapletters[1]{\wr@pletters#1\@nil}
\def\wr@pletters#1#2\@nil{#1\allowbreak\if&#2&\else\wr@pletters#2\@nil\fi}

\usepackage{fixme}
\fxsetup{
    status=draft,
    author=,
    layout=margin,
    theme=color
}

\begin{document}
\title[$k$-AP free sets can have many smaller progressions]{Sets without $k$-term progressions can have many shorter progressions}

\author{Jacob Fox}
\address{Department of Mathematics, Stanford University, Stanford, CA 94305, USA}
\email{jacobfox@stanford.edu}

\author{Cosmin Pohoata}
\address{Department of Mathematics, California Institute of Technology, Pasadena, CA 91106, USA}
\email{apohoata@caltech.edu}

\begin{abstract}


Let $f_{s,k}(n)$ be the maximum possible number of $s$-term arithmetic progressions in a set of $n$ integers which contains no $k$-term arithmetic progression. For all fixed integers $k > s \geq 3$, we prove that $f_{s,k}(n)=n^{2-o(1)}$, which answers an old question of Erd\H{o}s. In fact, we prove upper and lower bounds for $f_{s,k}(n)$ which show that its growth is closely related to the bounds in Szemer\'edi's theorem.
\end{abstract}
\maketitle
\section{Introduction}
\label{intro}

\bigskip

Let $k \geq 3$ be an integer. In this paper, a $k$-term arithmetic progression of integers will denote as usual a set of the form $\{x,x+d, \dots, x+(k-1)d\}$. If $d \neq 0$, then we say that the progression is {\it{non-trivial}}. If a set $A$ does not contain any non-trivial $k$-term arithmetic progressions, we say that $A$ is \textit{$k$-AP free}. The study of $k$-AP free sets in the integers and other groups has been a central topic in additive combinatorics. Following the standard notation, we will denote by $r_k(n)$ the size of the largest $k$-AP free subset of $\left\{1,\ldots,n\right\}$. The seminal result on this topic is Szemer\'{e}di's Theorem \cite{S}, which states that sets of integers with positive density contain arbitrarily long arithmetic progressions, or using the notation above $r_k(n)=o(n)$. 

Since Szemer\'{e}di, the problem of finding better quantitative bounds for $r_{k}(n)$ has received a lot of attention, with impressive progress that led to many important tools, which in the meantime have become standard. For our application, we won't need the best bounds for each $k$, so we will limit ourselves to only mentioning Gowers' theorem \cite{Gowers1, Gowers2} that for each $k \geq 3$ there exists an absolute constant $c_{k} > 0$ such that
\begin{equation} \label{Gowers}
r_{k}(n) \ll \frac{n}{(\log \log n)^{c_{k}}}.
\end{equation}
Regarding lower bounds, Rankin \cite{Rankin} showed that there exists a constant $c_{k}' > 0$ such that
\begin{equation} \label{Rankin}
r_{k}(n) \gg \frac{n}{2^{c_{k}'(\log n)^{1/\lceil \log k \rceil}}}.
\end{equation}
Throughout the paper, all logarithms are base $2$ and the signs $\ll$ and $\gg$ are the usual Vinogradov symbols. 

Let $\mathcal{A}_{k}(n)$ be the set of $n$-term nonnegative integer sequences which contain no $k$-term arithmetic progression as a subsequence. Furthermore, let $f_{s}(A)$ denote the number of $s$-term arithmetic progressions in $A$, and finally let $f_{s,k}(n) = \max_{A \in \mathcal{A}_{k}(n)} f_{s}(A)$. In \cite[page 119]{Erdos}, Erd\H{o}s observed that
$$\frac{\log f_{3,4}(n)}{ \log n} > 1.4649$$
holds for infinitely many $n$ by constructing examples of sequences $A \in \mathcal{A}_{4}(3^{s})$ for which $f_{3}(A) = 3^{s-1}$. Furthermore, he noticed that for each $k > 3$ the limit $\lim_{n \to \infty} \log f_{3,k}(n)/\log n:=f_{3,k}$ exists, and asked whether or not $f_{3,k}$ is always less than $2$. In \cite{SA}, Simmons and Abbott improved on Erd\H{o}s' observation by showing that $f_{3,4}(n) \geq n^{1.623}$ holds infinitely often, and also proved that $f_{3,k} \to 2$ as $k$ goes to infinity. Nonetheless, in the regime when $k$ is fixed, there has been no further progress on understanding the limit $f_{3,k}$ as far as we are aware of. In this note, we settle Erd\H{o}s' question in the negative by proving the following more general result.

\smallskip

\begin{theorem} \label{main}
For all integers $k > s \geq 3$, we have
$$\lim_{n \to \infty} \frac{\log f_{s,k}(n)}{\log n} = 2.$$
\end{theorem}

\smallskip

In fact, we prove upper and lower bounds for $f_{s,k}(n)$ which show that its growth is closely related to the bounds in Szemer\'edi's theorem.

\smallskip

\begin{theorem} \label{UpLow}
There exist absolute positive constants $c$ and $C$ such that, for integers $k > s \geq 3$ and every sufficiently large integer $n$, we have 
$$ \left(\frac{c \cdot r_{k}(n)}{n}\right)^{2(s-2)} \cdot n^{2}\leq f_{s,k}(n) \leq \left(\frac{r_{k}(n)}{n}\right)^{C} \cdot n^{2}.$$
\end{theorem}

\smallskip

In light of the bounds on $r_{k}(n)/n$ provided by \eqref{Gowers} and \eqref{Rankin}, it is easy to check that Theorem \ref{main} follows from Theorem \ref{UpLow}; therefore, it suffices to prove the latter. We will do this already in Section 2. The proof of Theorem \ref{UpLow} will require a few ingredients from additive combinatorics, but we will state them in full as we will get to apply them, as they do not require much preparation.

\bigskip

\subsection*{Funding and acknowledgments}

The first author was supported by a Packard Fellowship and by NSF grant DMS-1855635. The second author would like to thank Oliver Roche-Newton for helpful conversations.

\bigskip

\section{Proof of Theorem \ref{UpLow}}

\bigskip

We first prove the desired upper bound on $f_{s,k}(n)$. For $s \geq 3$, we have $f_{s,k}(n) \leq f_{3,k}(n)$, so in order to prove the upper bound it suffices to show that
$$f_{3,k}(n) \leq \left(\frac{r_{k}(n)}{n}\right)^{C} n^{2}$$
holds for some absolute constant $C > 0$ and sufficiently large $n$. We will in fact show this claim for $C = 1/25$. Let $A \in \mathcal{A}_{k}(n)$ and let $pn^{2}$ denote the number of three-term arithmetic progressions in $A$, where $p$ is some positive real number (which is strictly less than $1$); i.e. $f_{3}(A)=pn^2$. 

To upper bound $p$, we will require the following variant of the Balog-Szemer\'edi-Gowers theorem (see \cite[Proposition 7.3, page 503]{Gowers1} or \cite[Section 5.1]{FS}).

\begin{theorem} \label{BSG}
If $A$ and $B$ are sets of $n$ integers and $G$ is a bipartite graph between $A$ and $B$ with $pn^2$ edges such that partial sumset $A +_{G} B$ has size at most $K|A|$, then there is a subset $A'$ of $A$ with $|A'| \geq pn/4$ and 
$$|A'-A'| \ll K^4p^{-5}n.$$
\end{theorem}

Here $A +_{G} B$ denotes as usual the sumset restricted to the edges coming from $G$, namely
$$A +_{G} B = \left\{a+b:\ a \in A, b \in B, (a,b) \in E(G)\right\}.$$
It is perhaps important to mention that Theorem \ref{BSG} is a somewhat nonstandard version of the Balog-Szemer\'edi-Gowers theorem, which outputs directly a large set $A' \subset A$ with small difference set, without applying any Ruzsa-type inequality. One can derive this version from the following lemma.

\begin{lemma} \label{paths}
If a bipartite graph $G=(A,B,E)$ with $|A|=|B|=n$ has $pn^2$ edges with $p \geq 30n^{-1/2}$, then there is a subset $A'$ of $A$ of size at least $pn/4$ such that every pair of vertices in $A'$ have at least $\Omega(p^5n^3)$ paths of length four connecting them. 
\end{lemma}

For the sake of completeness, we include a quick proof Lemma \ref{paths}. We apply Lemma 5.1 from \cite{FS} with $\epsilon = 1/10$ and $c=p$ to obtain a subset $U$ of $A$ with $|U| \geq pn/2$ such that at least a $.9$ fraction of pairs of vertices in $U$ have at least $p^2n/20$ common neighbors. Consider the auxiliary graph $F$ on $U$ where two vertices are adjacent if in the original graph they have less than $p^2n/20$ common neighbors. By construction, the average degree in $F$ is at most $.1(|U|-1)$, so there are at most $|U|/2$ vertices of degree at most twice the average degree, which is at most $.2(|U|-1)$. Let $A'$ be the $|U|/2$ vertices of minimum degree in $F$. Then for every two vertices $a,a'$ in $A'$, their number of common neighbors $a''$ in the complement of F is at least $|U|-2-.4(|U|-1) \geq |U|/2$. Then any choice of this common neighbor $a''$ in F can be used as the middle vertex of at least $(p^2n/20)(p^2n/20 - 1)$ paths of length four between $a$ and $a'$ in the original graph, giving a total of at least $\frac{|U|}{2} (p^2n/20)(p^2n/20-1) \geq p^5 n^3 / 2000$ paths of length three between $a$ and $a'$. 

Using Lemma \ref{paths}, one can then deduce Theorem \ref{BSG} in the usual way. First observe that we may assume $p\geq 30n^{-1/2}$ since otherwise Theorem \ref{BSG} is trivial taking $A'=A$. Applied to the graph from the setup of Theorem \ref{BSG}, Lemma \ref{paths} produces $A' \subset A$ of size at least $pn/4$ such that every pair of vertices in $A$ have at least $\Omega(p^5n^3)$ paths of length four connecting them. This set happens to also satisfy $|A'-A'| \ll K^4p^{-5}n$. Indeed, for each $a, a' \in A'$, consider a path of length four in $G$ between them, say $(a,b,a'',b',a')$. For $y:=a-a' \in A'-A'$, we can then write
$$a-a' = (a+b) - (a''+b) + (a''+b') - (a'+b') = x_{1}-x_{2}+x_{3}-x_{4},$$
where $x_{1}=a+b$, $x_{2}=a''+b$, $x_{3}=a''+b'$, and $x_{4}=a'+b'$ are all elements of $A +_{G} B$. Since for every $a, a' \in A'$ there are at least $\Omega(p^5n^{3})$ paths of length four between $a$ and $a'$, this means every $y \in A'-A'$ can be written as $x_{1}-x_{2}+x_{3}-x_{4}$ for at least $\Omega(p^5 n^3)$ quadruples $(x_{1},x_{2},x_{3},x_{4}) \in (A+_{G} B)^{4}$. However, $|A +_{G} B| \leq Kn$ holds by assumption, so there are at most $K^{4}n^{4}$ such quadruples. By the pigeonhole principle, it then follows that the number of distinct elements $y \in A'-A'$ is at most $O(K^4p^{-5}n)$, as claimed.

Returning to the task of deriving the upper bound in Theorem \ref{UpLow}, we apply Theorem \ref{BSG} to the graph $G$ where $A$ and $B$ are chosen to be two copies of our $k$-AP free $A$ and with an edge between $(a,b) \in A \times A$ if $a+b=2c$ for some $c \in A$. This graph has precisely $pn^{2}$ edges and we can apply Theorem \ref{BSG} to it with $K=1$ since
$$|A +_{G} A| = |\left\{2a:\ a \in A\right\}| = |A|.$$
This yields a subset $A' \subset A$ with $|A'| \geq p|A|/4$ and $|A'-A'| \ll p^{-5}n \ll p^{-6}|A'|$. At this point, we recall a version of the so-called Freiman-Ruzsa modelling lemma (see for instance \cite[Theorem 2.3.5, page 127]{Ruzsa}).

\begin{lemma} \label{modelling}
Let $S$ be a finite set of integers and let $r \geq 2$ be an arbitrary integer. Then, there is a set $S^{*} \subset S$ with $|S^{*}| \geq |S|/r^{2}$ which is Freiman $r$-isomorphic to a set of integers $T$ such that
$$T \subset \left\{1,2,\ldots,\left\lceil\frac{1}{r} \cdot |rS-rS|\right\rceil\right\}.$$
\end{lemma}

Here $rS-rS$ denotes the sumset $S+\ldots+S-S-\ldots-S$, where $S$ appears $2r$ times. For the reader's convenience, we also recall that for any two commutative groups $G_{1}$, $G_{2}$ two sets $S \subset G_{1}$ and $T \subset G_{2}$ are said to be Freiman $r$-isomorphic if there exists a one to one map $\phi : S \to T$ such that for every $x_{1},\ldots,x_{r},y_{1},\ldots,y_{r}$ in $S$ (not necessarily distinct) the equation
$$x_{1}+\ldots+x_{r} = y_{1} + \ldots + y_{r}$$
holds if and only if
$$\phi(x_{1})+\ldots+\phi(x_{r}) = \phi(y_{1}) + \ldots + \phi(y_{r}).$$
We combine Lemma \ref{modelling} with (a consequence of) the classical Pl$\ddot{\text{u}}$nnecke-Ruzsa inequality, for which a simple proof can be found in \cite{Petridis}.
\begin{lemma}
Let $S$ and $T$ be finite sets of reals such that $|S+T| \leq \alpha |S|$, and let $r, r'$ be positive integers. Then
$$|rT - r'T| \leq \alpha^{r+r'}|S|.$$
\end{lemma}
Indeed, if we apply this with $S=A'$, $T=-A'$, $r=r'=2$, and $\alpha = p^{-6}$, we have
$$|2A'-2A'| \leq p^{-24} |A'|  \leq p^{-24}n.$$
Therefore, by Lemma \ref{modelling}, there is a subset $A^{*} \subset A'$ with $|A^{*}| \gg pn$ which is Freiman $2$-isomorphic to a set of integers $\phi(A^{*})$ contained in the interval $\left\{1,\ldots,\left\lceil p^{-24}n\right\rceil\right\}$. In particular, since $\phi$ preserves $k$-term arithmetic progressions,
$$pn \ll |A^{*}| = |\phi(A^{*})| \leq r_{k}(\left\lceil p^{-24}n\right\rceil).$$
Lastly, recall that $r_{k}(n)$ is subadditive as a function of $n$, namely the inequality $r_k(n+n') \leq r_k(n)+r_k(n')$ holds for all positive integers $n,n'$. In particular, $r_{k}(\left\lceil p^{-24}n\right\rceil) \ll p^{-24}r_{k}(n)$, hence $pn \ll p^{-24}r_{k}(n)$, or equivalently $p^{25} \ll r_{k}(n)/n$. This means that $A$ contains at most $\left(r_{k}(n)/n\right)^{1/25} n^{2}$ three-term arithmetic progressions. This completes the proof of the upper bound.

We next prove the desired lower bound on $f_{s,k}(n)$ in Theorem \ref{UpLow}. We begin by revisiting some further simple properties of $r_{k}(n)$ as a function of $n$. In addition to being subadditive, we also recall that $r_{k}(n)$ is an increasing function, so $r_k(m) \leq r_k(n)$ if $m \leq n$. Together these imply that if $n \geq m$, we have $r_k(n) \leq \lceil \frac{n}{m}\rceil r_k(m) \leq \frac{2n}{m} r_k(m)$, so 
\begin{equation}\label{almostdecreasing}
\frac{r_k(n)}{2n} \leq  \frac{r_k(m)}{m}.
\end{equation}

For all positive integers $m$ and $n$, we have 
\begin{equation} \label{supermult}
r_{k}(2mn) \geq r_{k}(m)r_{k}(n).
\end{equation}
Indeed, if $U$ is a subset of $\left\{1,\ldots,m\right\}$ without a $k$-term arithmetic progression and $V$ is a subset of $\left\{1,\ldots, n\right\}$ without a $k$-term arithmetic progression, then the set 
$$W=\left\{2u(n-1)+v:\ u \in U, v \in V\right\}$$
is a $k$-AP free subset of $\left\{1,\ldots,2mn\right\}$ of size $|U||V|$, so \eqref{supermult} follows. 

In particular, if $n \geq N^{1/2}$, letting $m=\lfloor \frac{N}{2n} \rfloor$, we have 
$$r_k(N) \geq r_k(2mn)\geq r_k(n)r_k(m) \geq r_k(n)\frac{m}{2n}r_k(n) \geq \frac{N}{8}\left(\frac{r_k(n)}{n}\right)^2,$$
where the first inequality follows from $r_k(n)$ being an increasing function, the second inequality is by (\ref{supermult}), the third inequality is by (\ref{almostdecreasing}) using $n \geq m$, and finally the fourth inequality is by substituting in $n \leq 4mN$. It thus follows that 
\begin{equation}\label{cor}
\frac{r_k(N)}{N} \geq \frac{1}{8}\left(\frac{r_k(n)}{n}\right)^{2}.
\end{equation}

Let $N=N_{n,k,s}$ be the least positive integer such that $r_{k}(N) = \lfloor n/s \rfloor$. Such an $N$ exists since, for every $m$, $r_k(m+1) =r_k(m)$ or $r_k(m)+1$ and $\lim_{m \to \infty} r_k(m)=\infty$. We will show that for $k > s \geq 3$ and $n$ sufficiently large in terms of $k$, we have 
\begin{equation} \label{lower}
f_{s,k}(n) \geq \left( \frac{n}{300sN} \right)^{s-2} n^2.
\end{equation}

For $n$ sufficiently large in terms of $k$, we have $n \geq N^{1/2}$ holds (for instance by \eqref{Rankin}), so \eqref{cor} implies that 
$n/N \geq s \cdot r_{k}(N)/N \geq s \cdot (1/8) \cdot (r_{k}(n)/n)^{2}$, and hence the lower bound from Theorem \ref{UpLow} follows from \eqref{lower}. We next prove \eqref{lower} using a probabilistic construction of a $k$-AP free set $A$ of $n$ integers with many $s$-term arithmetic progressions. 

For each $1 \leq i \leq s$, let $d_{i}$ be an integer chosen uniformly and independently at random from the set $\{1,\ldots,2N\}$. 
Let $S \subset \left\{1,\ldots,N\right\}$ be a $k$-AP free set of cardinality $r_{k}(N)=\lfloor n/s \rfloor$, and $S_{i}$ denote the translate $\left\{x + 6(i-1)N-1+d_{i}:\ x \in S\right\}$, i.e. $S_{i} := S + \left\{6(i-1)N-1+d_{i}\right\}$. 
 

Finally, let us consider the set $A\subset \left\{1,\ldots,6sN\right\}$ defined by
$$A := \bigcup_{i=1}^{s} S_{i}.$$
We first check that such a (random) set must be $k$-AP free. Indeed, the sets $S_{1},\ldots,S_s$ are pairwise disjoint since, for each $1 \leq i \leq s$, we have 
$$S_{i} \subset \left\{6(i-1)N+1,\ldots,6(i-1)N+3N-1\right\}.$$
Furthermore, these sets are spaced out so that if an arithmetic progression contains an element from $S_i$ and an element of $S_j$ with $i \not =j$, then its common difference is at least $3N+2$, in which case the arithmetic progression cannot contain two elements in the same $S_i$. In particular, every arithmetic progressions in $A$ of length longer than $s$ must be a subset one of the $S_i$, and hence $A$ is $k$-AP free.  Finally, $|A|=s|S|=s\lfloor \frac{n}{s} \rfloor \leq n$, so $A$ is indeed in $\mathcal{A}_{k}(n)$, or it can be artificially augmented to a set in $\mathcal{A}_{k}(n)$ by adding some elements that do not create $k$-term arithmetic progressions.


We next lower bound the expected number of $s$-term arithmetic progressions in $A$. The number of $s$-term arithmetic progressions $a,a+D,\ldots,a+(s-1)D$ with $a+(i-1)D \in \{6(i-1)N+N+1,\ldots,6(i-1)N+2N\}$ for $1 \leq i \leq s$ is the same as the number of $s$-term arithmetic progressions in $\{1,\ldots,N\}$ with any integer common difference, which is 
$$N+2\sum_{a=1}^{N-1}\left\lfloor \frac{N-a}{s}\right\rfloor \geq \frac{1}{s}{N \choose 2}.$$
For each such $s$-term arithmetic progression $a,a+D,\ldots,a+(s-1)D$ and for each sequence $(a_1,\ldots,a_s)$ of $s$ elements from $S$, there is a choice of $d_1,\ldots,d_s \in \{1,\ldots,2N\}$ such that $a_i+6(i-1)N-1+d_i=a+(i-1)D$ for $1 \leq i \leq s$. Hence, the expected number of $s$-term arithmetic progressions in $A$ is at least $$\frac{1}{s}{N \choose 2}|S|^s(2N)^{-s} \geq \frac{1}{4s}N^2\left(\frac{\lfloor n/s \rfloor}{2N}\right)^s \geq \left(\frac{n}{300sN}\right)^{s-2} n^2.$$  
Thus, there must exist a choice of such an $A$ for which the number of $s$-term arithmetic progressions is at least this lower bound on the expected number, which completes the proof of (\ref{lower}) and hence Theorem \ref{UpLow}.

\end{document}